\documentclass[letterpaper, 10 pt, conference]{ieeeconf}  

\IEEEoverridecommandlockouts
\usepackage{amsmath,amssymb,amsfonts}

\usepackage{bm}

\let\labelindent\relax
\usepackage{enumitem}

\usepackage{amsthm}

\theoremstyle{definition}
\newtheorem{definition}{Definition}[section]

\newtheorem{theorem}{Theorem}[section]

\newtheorem{remark}{Remark}[section]

\newtheorem{lemma}{Lemma}[section]

\newtheorem{proposition}{Proposition}[section]

\newtheorem{assumption}{Assumption}[section]

\usepackage{mathtools}

\usepackage{multicol}
\usepackage{multirow}
\usepackage{makecell}

\usepackage{pifont}
\newcommand{\tp}{^{\top}}

\newcommand{\iv}{^{-1}}
\newcommand{\tfin}{t_{\mathrm{f}}}

\newcommand{\norm}[1]{\left\lVert#1\right\rVert}

\makeatletter
\let\NAT@parse\undefined
\makeatother
\usepackage[noadjust]{cite}
\usepackage{xcolor}
\usepackage{tikz}
\usepackage{graphicx}
\newcommand{\dashedbox}[1]{\tikz[baseline=(X.base)]\node[draw,dashed,inner sep=2pt](X){$#1$};}
\usepackage[hidelinks]{hyperref}
\usepackage{cleveref} 
\crefname{assumption}{Assumption}{Assumptions}
\Crefname{assumption}{Assumption}{Assumptions}

\usepackage{cuted}
\usepackage{etoolbox}
\AfterEndEnvironment{strip}{\leavevmode}

\usepackage{dsfont}
\usepackage{nicefrac}
\usepackage{verbatim}

\usepackage{subcaption}
\usepackage{booktabs}
\usepackage{nth}
\usepackage{xfrac}

\newcommand{\bx}{\mathbf{x}}
\newcommand{\bu}{\mathbf{u}}
\newcommand{\bz}{\mathbf{z}}
\newcommand{\bv}{\mathbf{v}}
\newcommand{\bw}{\mathbf{w}}
\newcommand{\bff}{\mathbf{f}}
\newcommand{\ba}{\mathbf{a}}
\newcommand{\bc}{\mathbf{c}}
\newcommand{\bd}{\mathbf{d}}
\newcommand{\by}{\mathbf{y}}
\newcommand{\balpha}{\boldsymbol{\alpha}}
\newcommand{\bgamma}{\boldsymbol{\gamma}}
\newcommand{\bfeta}{\boldsymbol{\eta}}
\newcommand{\bzeta}{\boldsymbol{\zeta}}
\newcommand{\bxi}{\boldsymbol{\xi}}
\newcommand{\blambda}{\boldsymbol{\lambda}}
\newcommand{\bvarphi}{\boldsymbol{\varphi}}

\newcommand{\bGamma}{\boldsymbol{\Gamma}}
\newcommand{\mA}{\mathbf{A}}
\newcommand{\mB}{\mathbf{B}}
\newcommand{\mR}{\mathbf{R}}
\newcommand{\mG}{\mathbf{G}}
\newcommand{\mK}{\mathbf{K}}
\newcommand{\mW}{\mathbf{W}}
\newcommand{\mQ}{\mathbf{Q}}
\newcommand{\mF}{\mathbf{F}}
\newcommand{\mH}{\mathbf{H}}
\newcommand{\mI}{\mathbf{I}}
\newcommand{\mT}{\mathbf{T}}
\newcommand{\mL}{\mathbf{L}}
\newcommand{\cD}{\mathcal{D}}
\newcommand{\cC}{\mathcal{C}}
\newcommand{\cH}{\mathcal{H}}
\newcommand{\cN}{\mathcal{N}}
\newcommand{\cL}{\mathcal{L}}
\newcommand{\cP}{\mathcal{P}}
\newcommand{\cS}{\mathcal{S}}
\newcommand{\cV}{\mathcal{V}}
\newcommand{\cZ}{\mathcal{Z}}
\newcommand{\cE}{\mathcal{E}}

\newcommand{\cX}{\mathcal{X}}

\begin{document}

\title{Exact Finite-Dimensional Reduction of Waypoint-Constrained Trajectory Optimization for Feedback-Linearizable Systems \\
\thanks{This project has received funding from the European Union’s 2020
research and innovation programme under the Marie Sklodowska-
Curie grant agreement No. 953348 ELO-X. It was also supported by the Swiss National Science Foundation (SNSF) under Grant 200020\_219356 and by the European Research Council (ERC)
through European Union’s Horizon 2020 (H2020), under Grant Agreement 101020573 FunLearn. Corresponding author: {\tt\footnotesize shaohui.yang@epfl.ch}}
}

\author{Shaohui Yang$^{1,*}$, Vincent Guillemet$^{2,*}$, Bassam El Rawas$^{2}$, Michael Unser$^{2}$, and Colin N. Jones$^{1}$
\thanks{$^{*}$These authors contributed equally.}
\thanks{$^{1}$Automatic Control Laboratory, EPFL, Lausanne, Switzerland.}
\thanks{$^{2}$Biomedical Imaging Group, EPFL, Lausanne, Switzerland.}
}

\maketitle
\thispagestyle{plain}
\pagestyle{plain}

\begin{abstract}
Waypoint-constrained trajectory optimization (WCTO) seeks a minimum-effort
input that drives a system through prescribed states at finitely many
instants. For linear dynamics, its minimizers are splines with finitely many
coefficients, so the infinite-dimensional problem reduces exactly to a
finite-dimensional one. We extend this reduction to full-state feedback
linearizable systems with waypoint constraints imposed on the physical state.
In the linearizing coordinates, these constraints become manifolds that are
meaningful where the linearization is valid. We show that every local
minimizer is a finite kernel expansion, give a sufficient condition for the
existence of a minimizer, and certify the solutions of the finite-dimensional
program a posteriori. Experiments on a single-link arm confirm that the
reduction is exact. They also show that the linear model of the linearized
system is a design parameter. A proper choice can keep the solution where the
linearization is valid and can be used to modify the meaning of the cost being minimized on the linear space.
\end{abstract}

\section{Introduction}\label{sec:intro}

Many tasks in robotics and control specify where a system must be only at a few
instants, rather than along an entire trajectory. A quadrotor must
pass through a sequence of gates, a manipulator must reach grasping and release
poses at given times, or a vehicle must meet rendezvous conditions, while the
motion in between is left free. Among all motions that meet these requirements,
one naturally seeks the one with the least control effort. We refer to this problem as
waypoint-constrained trajectory optimization (WCTO), stated generically as
\begin{equation}\label{eq:wcto-text}
\begin{aligned}
    \min_{\bu(\cdot)}\ \ & \int_0^{\tfin}\bu(t)\tp \mR\, \bu(t)\,\mathrm{d}t\\
    \text{s.t.}\ \ & \dot \bx(t)=\bff\big(\bx(t),\bu(t)\big),\quad t\in[0,\tfin],\\
    & \bx(t_i)\in\cX_i,\quad i=1,\dots,N,
\end{aligned}
\end{equation}
where $\bx$ denotes the state, $\bu$ denotes the input, $t_i\in(0,\tfin]$ and $\tfin$
are given time instants, and $\cX_i$ are given constraint manifolds. Two
features set WCTO apart from generic optimal control: 1) the constraints act
only at finitely many prescribed instants; and 2) the decision variable is a
function of continuous time. Hence WCTO is an infinite-dimensional problem
with finitely many constraints.

For linear dynamics, the minimizers of WCTO are splines \cite{de1966splines, schultz1967splines, unser2005cardinal} with finitely many
coefficients \cite{zhang1997splines,egerstedt2010control,aubin2021linearly}.
This exact finite-dimensional reduction is well established. For nonlinear
systems, roboticists reach linear spline fitting through differential flatness, with waypoints
imposed on the flat outputs and their derivatives, which are linearizing
coordinates. However, most tasks prescribe waypoints in the physical states. In the
linearizing coordinates, such a waypoint becomes a manifold.
To the best of the authors' knowledge, no previous work has solved a problem of this type, where waypoint constraints are imposed directly on the physical states. 

In this paper, we extend the exact reduction to full-state feedback
linearizable systems, which are equivalent to linear systems under a
nonlinear state transformation and a feedback law. Our contributions are
threefold: 1) an exact finite-dimensional reduction of WCTO for these
systems, where the finitely many waypoint constraints are imposed on the
nonlinear physical coordinates rather than on the linear coordinates;
2) numerical experiments that validate the exactness of the reduction; and
3) numerical experiments showing that the system matrices $(\mA,\mB)$ of the
equivalent linear system offer degrees of freedom beyond the default
chain-of-integrator form. On a single-link arm, a choice of $(\mA,\mB)$ guided
by domain knowledge keeps the solution where feedback linearization is valid and better approximates the physical cost. Such a choice depends on the problem. 

In Section~\ref{sec:prelim}, we review feedback linearization and WCTO
for linear systems, together with its use in robotics.
In Section~\ref{sec:main}, we formulate WCTO for feedback-linearizable systems and
prove the exact finite-dimensional reduction. In Section~\ref{sec:num}, we
validate the reduction and study the choice of $(\mA,\mB)$ on a single-link
arm. Section~\ref{sec:conclusion} concludes the paper.

\emph{Notation.} Plain letters denote scalars. Bold lowercase letters denote
vectors and vector-valued functions. 
Bold uppercase letters denote
matrices and matrix-valued functions. Calligraphic letters denote sets and
function spaces. Two exceptions keep the letters customary in control. The
diffeomorphism $\mT$ of feedback linearization is vector-valued but
uppercase, and the matrix-valued function $\bgamma$ in \eqref{eq:fl-form} is
lowercase. 
Depending on context, $\norm{\cdot}$ denotes the $\cL^2$ norm of a vector-valued function or the $\ell_2$ norm of a finite-dimensional vector or matrix.

\section{Preliminaries and related works}\label{sec:prelim}

\subsection{Full-State Feedback Linearization}\label{sec:prelim-fl}

Consider the input-affine nonlinear system
\begin{equation}\label{eq:nl-sys}
    \dot{\bx} = \bff(\bx) + \mG(\bx)\,\bu,
\end{equation}
with state $\bx\in \cD$ and input $\bu\in\mathbb{R}^m$, where $\cD\subseteq\mathbb{R}^n$
and $\bff:\cD\to\mathbb{R}^n$, $\mG:\cD\to\mathbb{R}^{n\times m}$ are
sufficiently smooth.

\begin{definition}[{\cite[Def.~13.1]{khalil2002nonlinear}}]\label{def:fl}
The system \eqref{eq:nl-sys} is \emph{feedback linearizable} (FL), or
\emph{input-state linearizable}, if there exists a diffeomorphism
$\mT:\cD\to\mathbb{R}^n$ such that $\cD_z\coloneqq \mT(\cD)$ contains the origin and the
change of variables $\bz=\mT(\bx)$ transforms \eqref{eq:nl-sys} into
\begin{equation}\label{eq:fl-form}
    \dot{\bz} = \mA\bz + \mB\bgamma(\bx)\big[\bu-\balpha(\bx)\big],
\end{equation}
with $(\mA,\mB)\in\mathbb{R}^{n\times n}\times\mathbb{R}^{n\times m}$ controllable,
$\balpha:\cD\to\mathbb{R}^m$, and $\bgamma:\cD\to\mathbb{R}^{m\times m}$ nonsingular
for all $\bx\in \cD$.
\end{definition}

The feedback $\bu=\balpha(\bx)+\bgamma(\bx)\iv \bv$ yields the linear system
\begin{equation}\label{eq:lin-sys}
    \dot{\bz} = \mA\bz + \mB\bv, \qquad \bz\in \cD_z,
\end{equation}
with virtual input $\bv\in\mathbb{R}^m$. We call $(\bx, \bu)$ and $(\bz, \bv)$ the \emph{physical}
and \emph{linearizing} coordinates, respectively.

Textbooks usually present feedback linearization with $(\mA,\mB)$ in the
chain-of-integrator form, the Brunovsk\'y form
\cite{brunovsky1970classification}. However, \Cref{def:fl} only requires
$(\mA,\mB)$ to be controllable. The following proposition makes this freedom
explicit. It underlies contribution 3).

\begin{proposition}[The linear model is a design parameter]\label{prop:fl-AB}
Let $\mQ\in\mathbb{R}^{n\times n}$ and $\mH\in\mathbb{R}^{m\times m}$ be
nonsingular and $\mF\in\mathbb{R}^{m\times n}$. If
$(\mT,\mA,\mB,\balpha,\bgamma)$ satisfies \eqref{eq:fl-form}, then so does
$\big(\mQ\mT,\,\mQ(\mA+\mB\mF)\mQ\iv,\,\mQ\mB\mH,\,\balpha+\bgamma\iv \mF\mT,\,
\mH\iv\bgamma\big)$ in the coordinates $\bar\bz=\mQ\bz$ and $\bar\bv$ with
$\bv=\mF\bz+\mH\bar\bv$.
\end{proposition}

Hence, $(\mA,\mB)$ is a design parameter as long as
controllability is maintained.

\begin{remark}[No zero dynamics]\label{rem:fl-full}
\Cref{def:fl} linearizes all $n$ states. This corresponds to a relative
degree $n$ with respect to a suitable output and leaves no zero dynamics.
Hence, a waypoint on any part of the physical state $\bx$ has a counterpart in
the linearizing coordinates.
\end{remark}

\begin{remark}[Locality]\label{rem:fl-local}
Feedback linearization holds on $\cD \subset \mathbb{R}^n$ rather than the entire space. A trajectory in the linearizing coordinates defines a physical
trajectory only if it stays in $\cD_z$.
\end{remark}

\subsection{Waypoint-Constrained Trajectory Optimization}\label{sec:prelim-lin}

We start by describing WCTO \eqref{eq:wcto-text} with linear dynamics and constraints, where
the theory is well established. We first review the literature, then restate
the solution with tools from functional analysis, and finally review its use
for nonlinear systems in robotics.

\subsubsection{Linear WCTO and Its Literature}\label{sec:prelim-lit}

Given a final time $\tfin>0$, a controllable pair $(\mA,\mB)$ with $\mB$ of full
column rank, a weight matrix $\mR\succ0$, an initial state $\bz_0\in\mathbb{R}^n$,
waypoint data $(t_i,\bc_i,y_i)\in(0,\tfin]\times\mathbb{R}^n\times\mathbb{R}$,
$i=1,\dots,N$, and loss functions $E_i:\mathbb{R}\times\mathbb{R}\to
\mathbb{R}\cup\{+\infty\}$ that are lower semicontinuous and bounded below in
their second argument, linear WCTO reads
\begin{equation}\label{eq:plin}
\begin{aligned}
    \min_{\bv}\ \ &
    \int_0^{\tfin} \bv(t)\tp \mR\,\bv(t)\,\mathrm{d}t
    +\sum_{i=1}^N E_i\big(y_i,\,\bc_i\tp \bz(t_i)\big)\\
    \text{s.t.}\ \ & \dot{\bz}=\mA\bz+\mB\bv,\quad \bz(0)=\bz_0.
\end{aligned}
\end{equation}
Problem \eqref{eq:plin} enforces the waypoints of \eqref{eq:wcto-text} through
loss functions rather than constraint manifolds. For a set $\cC$, the indicator
$\iota_{\cC}$ equals $0$ on $\cC$ and $+\infty$ outside, so the loss
$E_i(y_i,\cdot)=\iota_{\{y_i\}}(\cdot)$ enforces $\bc_i\tp\bz(t_i)=y_i$. The
two formulations are used interchangeably throughout the paper.

For single-input systems, Zhang \emph{et al.} \cite{zhang1997splines} solved
\eqref{eq:plin} with Hilbert space techniques. With indicator losses, the minimizers
are interpolating splines \cite{zhang1997splines}. With quadratic penalties,
they are smoothing splines, and further constraints such as monotonicity
lead to quadratic programs (QPs) \cite{egerstedt2010control}.
Aubin-Frankowski \cite{aubin2021linearly} solved the multi-input,
time-varying case. He identified the trajectories of the linear system as a
vector-valued reproducing kernel Hilbert space (vRKHS), constructed its
kernel explicitly, and derived the minimizers from a representer theorem.
Continuous-time affine state constraints of the linear-quadratic regulator
(LQR) are enforced through tightenings at finitely many
times, which yields a second-order cone program (SOCP). Aubin-Frankowski thus
brought vRKHS to linear WCTO. Spline fitting for vector-valued functions dates back to the vector-valued $Lg$-splines of Sidhu and Weinert,
which cover interpolating \cite{sidhu1979vector} and smoothing
\cite{sidhu1984vector} splines. Under
generalized total-variation regularization, an $\cL^1$-type penalty \cite{unser2017splines}, the
counterpart of \eqref{eq:plin} has adaptive vector-valued splines in Banach
space as solutions \cite{guillemet2025adaptive}. 

\subsubsection{Solution via a Representer Theorem}\label{sec:prelim-math}

We restate the solution from \cite{aubin2021linearly} for time-invariant
dynamics. It views \eqref{eq:plin} as a minimum-norm problem over state
trajectories with finitely many evaluation losses. We first
construct the space of these trajectories.

We write the dynamics as $\mL_A\bz=\mB\bv$ with the differential operator
$\mL_A\coloneqq\frac{\mathrm{d}}{\mathrm{d}t}\mI_n-\mA$. Its null space
$\cN_A\coloneqq\{\bz:\ \mL_A\bz=\mathbf{0}\}=\{e^{\mA\cdot}\bzeta:\
\bzeta\in\mathbb{R}^n\}$ consists of the free responses.
Let $\mB^\dagger\coloneqq(\mB\tp \mB)\iv \mB\tp$ be the left pseudo-inverse of $\mB$,
and equip the square-integrable function space $\cL^2([0,\tfin];\mathbb{R}^m)$ with the inner product
$\langle \bv_1,\bv_2\rangle_R\coloneqq\int_0^{\tfin}\bv_1(t)\tp \mR\bv_2(t)\,\mathrm{d}t$
and the norm $\norm{\cdot}_R$. The projected inputs
$\cL^2_{R|B}\coloneqq\{\mB\bv:\ \bv\in \cL^2([0,\tfin];\mathbb{R}^m)\}$ form a Hilbert
space with the inner product
$\langle\cdot,\cdot\rangle_{R|B}\coloneqq\langle \mB^\dagger\cdot,\mB^\dagger\cdot\rangle_R$.
The state trajectories with finite cost form the space
$\cH_{A,B}\coloneqq\{\bz\ \text{absolutely continuous}:\ \mL_A\bz\in \cL^2_{R|B}\}$.
Each $\bz\in\cH_{A,B}$ determines the unique input
$\bv_{\bz}\coloneqq \mB^\dagger\mL_A\bz$. By the variation-of-constants formula,
$\bz(t)=e^{\mA t}\bz(0)+\int_0^{t}e^{\mA(t-r)}\mB\bv_{\bz}(r)\,\mathrm{d}r$. Hence
$\cH_{A,B}=\cN_A\oplus\cH^0_{A,B}$, where
$\cH^0_{A,B}\coloneqq\{\bz\in\cH_{A,B}:\bz(0)=\mathbf{0}\}$ collects the
responses from rest. The inner product
$\langle \bz_1,\bz_2\rangle\coloneqq \bz_1(0)\tp \bz_2(0)
+\langle\mL_A\bz_1,\mL_A\bz_2\rangle_{R|B}$ makes
$\cH_{A,B}$ a Hilbert space and the above sum orthogonal. Its second
term equals $\langle \bv_{\bz_1},\bv_{\bz_2}\rangle_R$.

\begin{definition}\label{def:vrkhs}
A Hilbert space $(\cH,\langle\cdot,\cdot\rangle)$ of
$\mathbb{R}^n$-valued functions on $[0,\tfin]$ is a \emph{vector-valued
Reproducing Kernel Hilbert Space} (vRKHS) if there is a matrix-valued kernel
$\mK:[0,\tfin]^2\to\mathbb{R}^{n\times n}$ such that, for all $t\in[0,\tfin]$ and
$\bfeta\in\mathbb{R}^n$, $\mK(\cdot,t)\bfeta\in\cH$ and the reproducing property
$\bfeta\tp \bw(t)=\langle \mK(\cdot,t)\bfeta,\bw\rangle$ holds for every $\bw\in\cH$.
\end{definition}

By the Riesz representation theorem, $\cH$ is a vRKHS if and only if
every evaluation $\bw\mapsto\bfeta\tp \bw(t)$ is continuous, and $\mK$ is then
unique. For $\bw\in\cH^0_{A,B}$, the variation-of-constants formula gives
$\bfeta\tp \bw(s)=\langle \mathds{1}_{[0,s]}\mR\iv \mB\tp
e^{\mA\tp(s-\cdot)}\bfeta,\,\bv_{\bw}\rangle_R$, which is continuous by the
Cauchy--Schwarz inequality. Its representer is the response from rest to
this input, which vanishes after $s$. Likewise, $\bfeta\tp \bz(s)=(e^{\mA\tp
s}\bfeta)\tp \bz(0)$ for $\bz\in\cN_A$. Thus $\mK=\mK_0+\mK_1$, where
    $\mK_0(t,s)\coloneqq e^{\mA t}e^{\mA\tp s}$
reproduces $\cN_A$ and
\begin{equation}\label{eq:K1}
    \mK_1(t,s)\coloneqq\int_0^{\min(t,s)}e^{\mA(t-r)}\mB\mR\iv \mB\tp
    e^{\mA\tp(s-r)}\,\mathrm{d}r
\end{equation}
reproduces $\cH^0_{A,B}$. $\mW(t)\coloneqq \mK_1(t,t)$ is the $\mR$-weighted controllability
Gramian on $[0,t]$.

Since $\bz_0$ is fixed, we write $\bz=e^{\mA\cdot}\bz_0+\bw$ with
$\bw\in\cH^0_{A,B}$ and $\norm{\bw}=\norm{\bv_{\bz}}_R$. Problem~\eqref{eq:plin}
is then equivalent to
\begin{equation}\label{eq:plin-w}
    \min_{\bw\in\cH^0_{A,B}}\ \ \norm{\bw}^2
    +\sum_{i=1}^NE_i\big(y_i,\hat y_i+\bc_i\tp \bw(t_i)\big),
\end{equation}
where $\hat y_i\coloneqq \bc_i\tp e^{\mA t_i}\bz_0$ is the waypoint value attained by
the free response. The dynamics now enter only through the constraint
$\bw\in\cH^0_{A,B}$, and the losses depend on finitely many evaluations
of $\bw$. Representer theorems \cite{wahba1990spline, scholkopf2001generalized} describe the solutions of optimization problems of the form \eqref{eq:plin-w}. Theorem 1 in \cite{aubin2021linearly}
exploits the specific structure of \eqref{eq:plin-w}.

\begin{theorem}[Representer theorem]\label{thm:rep-lin}
Let $\bGamma_c\in\mathbb{R}^{N\times N}$ be the Gram matrix with entries
$(\bGamma_c)_{ij}\coloneqq \bc_i\tp \mK_1(t_i,t_j)\bc_j$. If \eqref{eq:plin}
has a finite value, it has a minimizer, and every minimizer generates the
state trajectory
\begin{equation}\label{eq:plin-sol}
    \bz^\star(t)=e^{\mA t}\bz_0+\sum_{i=1}^N\lambda^\star_i \mK_1(t,t_i)\,\bc_i,
\end{equation}
where $\blambda^\star$ is a minimizer of
\begin{equation}\label{eq:plin-fin}
    \min_{\blambda\in\mathbb{R}^N}\ \blambda\tp\bGamma_c\blambda
    +\sum_{i=1}^N E_i\big(y_i,\hat y_i+(\bGamma_c\blambda)_i\big).
\end{equation}
If every $E_i(y_i,\cdot)$ is convex, then \eqref{eq:plin-fin}
is convex, and every local minimizer of \eqref{eq:plin-fin} generates a
global minimizer of \eqref{eq:plin}.
\end{theorem}

With the indicator loss $E_i(y_i,\sigma)=\iota_{\{y_i\}}(\sigma)$, \eqref{eq:plin-fin} reduces to the linear system
$\bGamma_c\blambda=\by-\hat\by$.
With the quadratic
penalty $E_i(y,\sigma)=\mu(y-\sigma)^2$, $\mu>0$, \eqref{eq:plin-fin} is a convex QP.
The shape of the spline \eqref{eq:plin-sol} is set by the eigenvalues of $\mA$
and $-\mA$ \cite{zhang1997splines}. A nilpotent $\mA$ gives polynomial splines,
and eigenvalues with nonzero real or imaginary parts add exponential or
trigonometric factors.

\subsubsection{WCTO with Nonlinear Dynamics in Robotics}\label{sec:prelim-robotics}

In motion planning, WCTO is commonly applied to nonlinear systems through
differential flatness. The flat outputs and their derivatives form
linearizing coordinates $\bz$ as in \Cref{def:fl}, possibly after a dynamic
extension of the state. Since these coordinates are successive derivatives,
$(\mA,\mB)$ is in chain-of-integrator form, and the virtual input $\bv$ consists of
the highest derivatives of the flat outputs. Minimizing the integrated squared
$k$th derivative yields polynomial splines of degree $2k-1$. Minimum
acceleration gives cubic splines, minimum jerk quintic splines, minimum snap
septic splines, and so on. The vectors $\bc_i$ select the flat outputs and
derivatives at $t_i$, so the waypoint constraints are affine in
$\bz$.

For quadrotors, minimum-snap trajectory generation
\cite{mellinger2011minimum} and MINCO \cite{wang2022geometrically} prescribe
the position, a flat output, at the waypoints and minimize the integrated
squared snap. Polynomial splines are optimal in the
infinite-dimensional setting of MINCO \cite[Thm.~2]{wang2022geometrically}. Industrial
manipulators connect joint-space waypoints with cubic splines
\cite{lin1983formulation}, minimum-jerk via-point trajectories reproduce human
reaching movements \cite{flash1985coordination}, and automated vehicles plan
lane changes with jerk-optimal quintic polynomials in a Frenet frame
\cite{werling2010optimal}. In all these works, the waypoints are affine in the
linearizing coordinates, so the linear theory of Section~\ref{sec:prelim-math} applies directly.

\section{Main Results}\label{sec:main}

In Section~\ref{sec:formulation}, we state WCTO for FL systems with waypoints prescribed
in physical coordinates and rewrite it in the linearizing coordinates. We then recast it as an equivalent loss-regularized problem in Section~\ref{sec:loss}
and prove a representer theorem.
The conditions for the existence of a global
minimizer are given in Section~\ref{sec:existence}. In Section~\ref{sec:reduction}, we show that every local minimizer is a
finite kernel expansion, reducing the problem exactly to a
finite-dimensional program.

\subsection{Problem Formulation}\label{sec:formulation}

We pose WCTO for the FL system \eqref{eq:nl-sys}. By \Cref{rem:fl-full},
waypoints on any part of the physical state $\bx$ then have a counterpart in
$\bz$. Given $\tfin>0$, $\mR\succ0$,
an initial state $\bx_0\in \cD$, and waypoint data
$(t_i,\bc_i,y_i,\diamond_i)\in(0,\tfin]\times(\mathbb{R}^n\setminus\{\mathbf{0}\})
\times\mathbb{R}\times\{=,\le\}$, $i=1,\dots,N$, with $t_1\le\dots\le t_N$, the WCTO problem for FL systems is
\begin{equation}\label{eq:wcto-fl}
\begin{aligned}
    \min_{\bu}\ \ & \int_0^{\tfin}\bv(t)\tp \mR\,\bv(t)\,\mathrm{d}t,\quad
    \bv\coloneqq\bgamma(\bx)\big[\bu-\balpha(\bx)\big],\\
    \text{s.t.}\ \ & \dot \bx=\bff(\bx)+\mG(\bx)\bu,\quad \bx(0)=\bx_0,\quad \bx(t)\in \cD,\\
    & \bc_i\tp \bx(t_i)\diamond_i y_i,\quad i=1,\dots,N.
\end{aligned}
\end{equation}
The cost penalizes the virtual input $\bv$ of the linearizing coordinates
rather than the physical input $\bu$. This is common practice in trajectory planning \cite{biagiotti2008trajectory} and optimal control \cite{beaver2024optimal}.
Its meaning depends on the choice of
$(\mA,\mB)$ in \Cref{prop:fl-AB}. In Section~\ref{sec:num-model}, we give a choice
for which $\bv$ approximates a scaled $\bu$ near an operating point.

Let $\bz_0\coloneqq\mT(\bx_0)$. By \Cref{def:fl}, $\bz=\mT(\bx)$ and
$\bu=\balpha(\bx)+\bgamma(\bx)\iv \bv$ map the feasible solutions of
\eqref{eq:wcto-fl} with finite cost bijectively onto the feasible solutions of
\begin{equation}\label{eq:P}\tag{P}
\begin{aligned}
    \min_{\bv}\ \ & \int_0^{\tfin}\bv(t)\tp \mR\,\bv(t)\,\mathrm{d}t\\
    \text{s.t.}\ \ & \dot \bz=\mA\bz+\mB\bv,\quad \bz(0)=\bz_0,\quad \bz(t)\in \cD_z,\\
    & \bc_i\tp \mT\iv\big(\bz(t_i)\big)\diamond_i y_i,\quad i=1,\dots,N. 
\end{aligned}
\end{equation}

Problem \eqref{eq:P} differs from \eqref{eq:plin} and from
all the robotics works of Section~\ref{sec:prelim-robotics} in one essential
respect. The affine waypoint constraints are imposed in the original nonlinear
coordinates $\bx$, not in the linearizing coordinates $\bz$. In $\bz$, the feasible
set of the $i$th waypoint,
\begin{equation}\label{eq:Ci}
    \cC_i\coloneqq\big\{\bzeta\in \cD_z:\ \bc_i\tp \mT\iv(\bzeta)\diamond_i y_i\big\},
\end{equation}
equals $\mT(\{\bx\in \cD:\ \bc_i\tp \bx\diamond_i y_i\})$ and is a subset of $\cD_z$.
Since $\mT\iv$ is undefined
outside $\cD_z$, a waypoint can only be prescribed where the linearization is
valid. The gradient
of $\bzeta\mapsto\bc_i\tp \mT\iv(\bzeta)$ never vanishes, since $\bc_i\neq\mathbf{0}$ and
$\mT\iv$ is a diffeomorphism. Hence $\cC_i$ is an embedded submanifold of
$\cD_z$, of dimension $n-1$ for an equality and of dimension $n$ with
boundary for an inequality. These manifolds are usually not flat, so
\eqref{eq:P} can be nonconvex. The path constraint
$\bz(t)\in \cD_z$ restates \Cref{rem:fl-local}, and we keep it throughout because
$\mT$ is rarely a global diffeomorphism.

\begin{remark}[Role of the selection vectors]\label{rem:selection}
For notational simplicity, the waypoint data $(t_i,\bc_i,y_i,\diamond_i)$ aim at scalar constraints. However, \eqref{eq:P} also supports vector constraints at $t_i$. 
If the corresponding equality constraint vectors $\bc_i$ span
$\mathbb{R}^n$, they prescribe the entire physical state $\bx(t_i)=\bar \bx_i$.
Then $\mT\iv(\bz(t_i))=\bar \bx_i$ is equivalent to the affine constraint
$\bz(t_i)=\mT(\bar \bx_i)$, and \eqref{eq:P} reduces to the linear WCTO
\eqref{eq:plin} with indicator losses. The problem becomes nontrivial when the vectors $\bc_i$
prescribe only part of the physical state, so that the waypoint set at $t_i$
is a manifold of positive dimension in $\bz$ rather than a point.
\end{remark}

\subsection{Loss Formulation and a Representer Theorem}\label{sec:loss}

We now rewrite \eqref{eq:P} as an equivalent loss-regularized problem.
Write $\bz=e^{\mA\cdot}\bz_0+\bw$ with $\bw\in\cH^0_{A,B}$ as in
Section~\ref{sec:prelim-math}, so that $\norm{\bw}=\norm{\bv_{\bz}}_R$, and set
\begin{equation}\label{eq:Ei}
    \ell_i(\bzeta)\coloneqq\iota_{\cC_i}\big(\bzeta+e^{\mA t_i}\bz_0\big)
\end{equation}
with the indicator of Section~\ref{sec:prelim-lit}.
The loss-regularized problem reads
\begin{equation}\label{eq:PE}\tag{P$_\ell$}
    \min_{\bw\in\cH^0_{A,B}}\ \norm{\bw}^2+\sum_{i=1}^N\ell_i\big(\bw(t_i)\big)
    \ \ \dashedbox{\text{s.t.}\ \ e^{\mA\cdot}\bz_0+\bw\in\cP},
\end{equation}
where the dashed box marks the path constraint with
$\cP\coloneqq\{\bz:\ \bz(t)\in \cD_z,\ \forall t\in[0,\tfin]\}$.

\begin{proposition}\label{prop:equiv}
The map $\bw\mapsto e^{\mA\cdot}\bz_0+\bw$ is a bijection from the feasible
points of \eqref{eq:PE} with finite cost onto the feasible points of
\eqref{eq:P}, and it preserves both the cost and the distance
$\norm{\bv_{\bz}-\bv_{\bz'}}_R$. Hence $\bw^\star$ is a local (resp.\ global)
minimizer of \eqref{eq:PE} if and only if $e^{\mA\cdot}\bz_0+\bw^\star$ is a
local (resp.\ global) minimizer of \eqref{eq:P}.
\end{proposition}

We therefore refer to \eqref{eq:P} and \eqref{eq:PE} interchangeably. Without
the path constraint, \eqref{eq:PE} can be tackled by the following theorem.
Its proof follows the standard argument \cite{aubin2021linearly}.

\begin{theorem}[Representer theorem with vector losses]\label{thm:rep-fl}
For $i=1,\dots,N$, let $\ell_i:\mathbb{R}^n\to\mathbb{R}\cup\{+\infty\}$ be lower
semicontinuous and bounded below, and let
$J_\ell(\bw)\coloneqq\norm{\bw}^2+\sum_{i=1}^N\ell_i\big(\bw(t_i)\big)$ for
$\bw\in\cH^0_{A,B}$. If $J_\ell$ is finite for some $\bw$, then $J_\ell$ has a
global minimizer, and every global minimizer $\bw^\star$ generates the state
trajectory
\begin{equation}\label{eq:rep-fl}
    \bz^\star(t)=e^{\mA t}\bz_0+\sum_{i=1}^N\mK_1(t,t_i)\,\ba_i^\star,\qquad
    \ba_i^\star\in\mathbb{R}^n.
\end{equation}
\end{theorem}

\begin{proof}
\emph{Existence.} Let $-\beta$ be a common lower bound of the $\ell_i$. A
minimizing sequence $(\bw_j)$ satisfies $\norm{\bw_j}^2\le J_\ell(\bw_j)+N\beta$, so it
is bounded and has a subsequence converging weakly to some $\bw^\star$. By the
reproducing property, $\bfeta\tp \bw_j(t)=\langle \mK_1(\cdot,t)\bfeta,\bw_j\rangle\to\bfeta\tp
\bw^\star(t)$ for all $t$ and $\bfeta$. Lower semicontinuity of the $\ell_i$ and weak
lower semicontinuity of the norm give $J_\ell(\bw^\star)\le\liminf_jJ_\ell(\bw_j)$.

\emph{Representer form.} Let $\cV$ be the span of
$\mK_1(\cdot,t_i)\bfeta$ over all $i$ and $\bfeta\in\mathbb{R}^n$, and
split a global minimizer as $\bw=\bw_{\cV}+\bw_\perp$ with $\bw_\perp\perp \cV$. The
reproducing property gives $\bfeta\tp \bw_\perp(t_i)=\langle
\mK_1(\cdot,t_i)\bfeta,\bw_\perp\rangle=0$ for all $\bfeta$, so $\bw(t_i)=\bw_{\cV}(t_i)$ and
$J_\ell(\bw)=J_\ell(\bw_{\cV})+\norm{\bw_\perp}^2$. Minimality gives
$J_\ell(\bw)\le J_\ell(\bw_{\cV})$, and $J_\ell(\bw)$ is finite. Hence
$\norm{\bw_\perp}^2\le0$, so $\bw_\perp=\mathbf{0}$ and $\bw\in\cV$.
\end{proof}

\Cref{thm:rep-fl} differs from \Cref{thm:rep-lin} in the losses. Each $\ell_i$
acts on the full evaluation $\bw(t_i)$ rather than on $\bc_i\tp \bw(t_i)$.
Hence each waypoint contributes a coefficient $\ba_i\in\mathbb{R}^n$ rather than
a scalar multiple of $\mK_1(\cdot,t_i)\bc_i$. 

\Cref{thm:rep-fl} requires lower semicontinuous losses, which here means
closed $\cC_i$. If $\mT$ is a global diffeomorphism, then $\cD_z=\mathbb{R}^n$, each
$\cC_i$ is closed, the path constraint is void, and \Cref{thm:rep-fl} yields
existence and the representer form at once. Two obstacles remain: 1) the set
$\cC_i$ is only closed relative to $\cD_z$, since its closure may contain boundary
points of $\cD_z$; and 2) $e^{\mA\cdot}\bz_0+\bw\in\cP$ requires infinitely many checks.

\subsection{Sufficient Condition for the Existence of a Minimizer}\label{sec:existence}

We overcome both obstacles by confining the trajectories of bounded
energy to a compact subset of $\cD_z$. The Gramian $\mW(t)=\mK_1(t,t)$ defined in
Section~\ref{sec:prelim-math} is positive definite for $t>0$ by
controllability. It is also nondecreasing in $t$ in the positive
semidefinite order, $\mW(t')\succeq \mW(t)$ for $t'\ge t$, since
$\mW(t)=\int_0^te^{\mA s}\mB\mR\iv \mB\tp e^{\mA\tp s}\,\mathrm{d}s$ has a positive
semidefinite integrand. For $\varepsilon\ge0$, define the ellipsoid
\begin{equation}\label{eq:ellipsoid}
    \cE_t(\varepsilon)\coloneqq\big\{e^{\mA t}\bz_0+\mW(t)^{1/2}\bxi:\
    \norm{\bxi}\le\varepsilon\big\}.
\end{equation}
The following lemma is the classical bound on the states reachable with
bounded input energy \cite{kalman1963controllability}, rederived with the
reproducing property.

\begin{lemma}\label{lem:tube}
Let $\bz\in\cH_{A,B}$ with $\bz(0)=\bz_0$. Then
$\bz(t)\in\cE_t(\norm{\bv_{\bz}}_R)$ for every $t\in[0,\tfin]$, and in
particular
\begin{equation}\label{eq:tube-ball}
    \norm{\bz(t)-e^{\mA t}\bz_0}\le\norm{\mW(t)}^{1/2}\,\norm{\bv_{\bz}}_R.
\end{equation}
\end{lemma}

\begin{proof}
Let $\bw\coloneqq \bz-e^{\mA\cdot}\bz_0\in\cH^0_{A,B}$, so that
$\norm{\bw}=\norm{\bv_{\bz}}_R$. The claims hold at $t=0$ since $\bw(0)=\mathbf{0}$, so let
$t>0$. By the reproducing
property and the Cauchy--Schwarz inequality,
$|\bfeta\tp \bw(t)|=|\langle \mK_1(\cdot,t)\bfeta,\bw\rangle|\le(\bfeta\tp \mW(t)\bfeta)^{1/2}\norm{\bw}$
for all $\bfeta$, because $\norm{\mK_1(\cdot,t)\bfeta}^2=\bfeta\tp \mW(t)\bfeta$. Taking
$\bfeta=\mW(t)\iv \bw(t)$ gives $\norm{\bxi}\le\norm{\bw}$ with
$\bxi\coloneqq \mW(t)^{-1/2}\bw(t)$, which is $\bz(t)\in\cE_t(\norm{\bv_{\bz}}_R)$.
Then $\norm{\bw(t)}^2=\bxi\tp \mW(t)\bxi\le\norm{\mW(t)}\norm{\bxi}^2$ gives
\eqref{eq:tube-ball}.
\end{proof}

Any feasible $\bar \bz$ bounds the cost of a minimizer by
$\norm{\bv_{\bar \bz}}_R^2$. By \Cref{lem:tube}, all trajectories within this
bound, including the tails of minimizing sequences, stay in the compact tube
of ellipsoids $\cE_t(\norm{\bv_{\bar \bz}}_R)$, $t\in[0,\tfin]$. If the tube
lies in $\cD_z$, it keeps a positive distance from the boundary of $\cD_z$,
where minimizers can fail to exist. The following assumption asks for exactly
this.

\begin{assumption}[Energy tube]\label{ass:tube}
There is a feasible $\bar \bz$ of \eqref{eq:P} such that, with
$\bar\varepsilon\coloneqq\norm{\bv_{\bar \bz}}_R$, the tube
$\cZ\coloneqq\bigcup_{t\in[0,\tfin]}\cE_t(\bar\varepsilon)$
satisfies $\cZ\subset \cD_z$.
\end{assumption}

\Cref{ass:tube} is hard to verify
beforehand, since it requires a feasible trajectory whose energy is small
enough for the whole tube to fit in $\cD_z$. The check simplifies if
$\mA\bz_0=\mathbf{0}$. The ellipsoids then share the center $\bz_0$ and are
nested because $\mW$ is nondecreasing, so $\cZ=\cE_{\tfin}(\bar\varepsilon)$.
Section~\ref{sec:reduction} shows how the strong assumption is bypassed in practice.

\begin{theorem}\label{thm:existence}
Under \Cref{ass:tube}, \eqref{eq:P} admits a global minimizer, and every
global minimizer $\bz^\star$ satisfies $\bz^\star(t)\in\cE_t(\bar\varepsilon)$
for all $t\in[0,\tfin]$.
\end{theorem}

\begin{proof}
The set $\cZ$ is the image of the compact set
$[0,\tfin]\times\{\norm{\bxi}\le1\}$ under the continuous map
$(t,\bxi)\mapsto e^{\mA t}\bz_0+\bar\varepsilon \mW(t)^{1/2}\bxi$, hence it is compact. Since
$\bc_i\tp \mT\iv$ is continuous on $\cZ\subset \cD_z$, the set
$\cC_i\cap\cZ$ is the preimage of the closed set $\{y_i\}$ or
$(-\infty,y_i]$ under a continuous map on a closed set, hence it is closed. The
indicator of a closed set is lower semicontinuous, so the auxiliary problem
\begin{equation}\label{eq:Paux}\tag{$\tilde{\mathrm{P}}$}
    \min_{\bw\in\cH^0_{A,B}}\ \norm{\bw}^2+\sum_{i=1}^N
    \iota_{\cC_i\cap\cZ}\big(\bw(t_i)+e^{\mA t_i}\bz_0\big),
\end{equation}
which is \eqref{eq:PE} with $\cC_i$ replaced by $\cC_i\cap\cZ$ and without
the path constraint, admits a global minimizer by \Cref{thm:rep-fl}. We
identify $\bw$ with $\bz=e^{\mA\cdot}\bz_0+\bw$ as in \Cref{prop:equiv} and compare
\eqref{eq:Paux} with \eqref{eq:P}. Let
$\cS\coloneqq\{\bz:\ \bz(0)=\bz_0,\ \norm{\bv_{\bz}}_R\le\bar\varepsilon\}$. By
\Cref{lem:tube}, every $\bz\in\cS$ satisfies
$\bz(t)\in\cE_t(\bar\varepsilon)\subset\cZ\subset \cD_z$ for all
$t$, so on $\cS$ the path constraint holds and $\bz(t_i)\in \cC_i$ if and
only if $\bz(t_i)\in \cC_i\cap\cZ$. Hence \eqref{eq:P} and
\eqref{eq:Paux} have the same feasible points and costs in $\cS$. Since
$\bar \bz$ is feasible for both with cost $\bar\varepsilon^2$, every minimizer of
either problem has cost at most $\bar\varepsilon^2$ and thus lies in
$\cS$. Thus \eqref{eq:P} and \eqref{eq:Paux} have the same minimizers,
so \eqref{eq:P} admits a global minimizer,
and every global minimizer $\bz^\star$ lies in $\cS$. \Cref{lem:tube}
then gives $\bz^\star(t)\in\cE_t(\bar\varepsilon)$ for all $t$.
\end{proof}

\subsection{Exact Finite-Dimensional Reduction}\label{sec:reduction}

\Cref{thm:existence} concerns global minimizers, but \eqref{eq:P} can be
nonconvex because the manifolds $\cC_i$ are not flat, and local solvers return local
minimizers. We now show that every local minimizer of \eqref{eq:P} is a finite
kernel expansion, so that \eqref{eq:P} reduces exactly to a finite-dimensional
program. Let $\ba\coloneqq(\ba_1,\dots,\ba_N)\in\mathbb{R}^{nN}$, let
$\bGamma\in\mathbb{R}^{nN\times nN}$ be the block Gram matrix with blocks
$\bGamma_{ij}\coloneqq \mK_1(t_i,t_j)$, let
$\bvarphi(\ba)\coloneqq\sum_j\mK_1(\cdot,t_j)\ba_j$ and
$\bzeta_i(\ba)\coloneqq e^{\mA t_i}\bz_0+\sum_j\mK_1(t_i,t_j)\ba_j$, and consider
\begin{equation}\label{eq:finite}
    \min_{\ba\in\mathbb{R}^{nN}}\ \ba\tp\bGamma\ba\quad
    \text{s.t.}\ \ \bzeta_i(\ba)\in \cC_i,\ \ i=1,\dots,N.
\end{equation}

\begin{theorem}\label{thm:reduction}
Let $\bz^\star$ be a local minimizer of \eqref{eq:P}. Then the following hold.
\begin{enumerate}[label=\arabic*),leftmargin=*]
\item $\bz^\star=e^{\mA\cdot}\bz_0+\bvarphi(\ba^\star)$ for some
$\ba^\star\in\mathbb{R}^{nN}$.
\item Every such $\ba^\star$ is a local minimizer of \eqref{eq:finite}.
\end{enumerate}
Conversely, if $\ba^\star$ is a local (resp.\ global) minimizer of
\eqref{eq:finite} and $\bz^\star\coloneqq e^{\mA\cdot}\bz_0+\bvarphi(\ba^\star)$ satisfies
$\bz^\star(t)\in \cD_z$ for all $t$, then $\bz^\star$ is a local (resp.\ global)
minimizer of \eqref{eq:P}.
\end{theorem}

\begin{proof}
\emph{Step 1 (margin).} The set $\bz^\star([0,\tfin])$ is compact and contained in
the open set $\cD_z$, so its distance $\rho$ to $\mathbb{R}^n\setminus \cD_z$ is
positive ($\rho=\infty$ if $\cD_z=\mathbb{R}^n$). Let $\bz$ satisfy $\bz(0)=\bz_0$.
Applying \eqref{eq:tube-ball} of \Cref{lem:tube} to $\bz-\bz^\star$, which starts at
the origin, and using that $\mW$ is nondecreasing gives
$\max_t\norm{\bz(t)-\bz^\star(t)}\le\kappa\norm{\bv_{\bz}-\bv_{\bz^\star}}_R$ with
$\kappa\coloneqq\norm{\mW(\tfin)}^{1/2}$. Hence every such
$\bz$ with $\norm{\bv_{\bz}-\bv_{\bz^\star}}_R<\rho/\kappa$ satisfies the
path constraint.

\emph{Step 2 (representer form).} Let $\cV\coloneqq\bvarphi(\mathbb{R}^{nN})$ and
split $\bw^\star\coloneqq \bz^\star-e^{\mA\cdot}\bz_0=\bw_{\cV}+\bw_\perp$ as in the proof of
\Cref{thm:rep-fl}. For $\nu\in[0,1]$, let
$\bz_\nu\coloneqq e^{\mA\cdot}\bz_0+\bw_{\cV}+\nu \bw_\perp$. Then $\bz_\nu(t_i)=\bz^\star(t_i)$, so the
waypoint constraints hold, and $\norm{\bv_{\bz_\nu}-\bv_{\bz^\star}}_R=(1-\nu)
\norm{\bw_\perp}$, so the path constraint holds for $\nu$ close to $1$ by Step~1.
If $\bw_\perp\neq\mathbf{0}$, the cost $\norm{\bw_{\cV}}^2+\nu^2\norm{\bw_\perp}^2$ of $\bz_\nu$ is
smaller than that of $\bz^\star$ for all $\nu<1$, contradicting local minimality.
Hence $\bw_\perp=\mathbf{0}$ and $\bw^\star=\bvarphi(\ba^\star)$ for some $\ba^\star$.

\emph{Step 3 (local minimality of $\ba^\star$).} Let $\ba$ be feasible for
\eqref{eq:finite}. Then $\bz_{\ba}\coloneqq e^{\mA\cdot}\bz_0+\bvarphi(\ba)$ satisfies
$\bz_{\ba}(t_i)=\bzeta_i(\ba)\in \cC_i$ and, by the reproducing property,
$\norm{\bv_{\bz_{\ba}}-\bv_{\bz^\star}}_R^2=(\ba-\ba^\star)\tp\bGamma(\ba-\ba^\star)$. For
$\norm{\ba-\ba^\star}$ small, Step~1 shows that $\bz_{\ba}$ is feasible for
\eqref{eq:P} and close to $\bz^\star$, so
$\ba\tp\bGamma\ba=\norm{\bv_{\bz_{\ba}}}_R^2\ge\norm{\bv_{\bz^\star}}_R^2
=\ba^{\star,\top}\bGamma\ba^\star$.

\emph{Step 4 (converse).} The trajectory $\bz^\star$ is feasible by assumption.
Let $\bz$ be feasible for \eqref{eq:P}, let $\bw\coloneqq \bz-e^{\mA\cdot}\bz_0$, and
let $\bw_{\cV}$ be its projection onto $\cV$. Since $\bvarphi$ is a linear map onto $\cV$
between finite-dimensional spaces, there are $\kappa_\varphi>0$ and $\ba$ with
$\bvarphi(\ba)=\bw_{\cV}$ and $\norm{\ba-\ba^\star}\le\kappa_\varphi\norm{\bw_{\cV}-\bvarphi(\ba^\star)}\le
\kappa_\varphi\norm{\bw-\bvarphi(\ba^\star)}$. As $\bzeta_i(\ba)=e^{\mA t_i}\bz_0+\bw_{\cV}(t_i)=\bz(t_i)\in \cC_i$, $\ba$
is feasible for \eqref{eq:finite}, and $\norm{\bv_{\bz}}_R^2=\norm{\bw}^2\ge
\norm{\bw_{\cV}}^2=\ba\tp\bGamma\ba$. If $\ba^\star$ is a global minimizer, the
right-hand side is at least $\ba^{\star\top}\bGamma\ba^\star=
\norm{\bv_{\bz^\star}}_R^2$. If $\ba^\star$ is a local minimizer, the same holds
whenever $\norm{\bv_{\bz}-\bv_{\bz^\star}}_R=\norm{\bw-\bvarphi(\ba^\star)}$ is small.
\end{proof}

Item 1) of \Cref{thm:reduction} extends the representer form \eqref{eq:rep-fl}
from global minimizers without the path constraint to local minimizers of
\eqref{eq:P}. \Cref{thm:reduction} states necessary conditions for local minimizers, whose
existence is guaranteed by \Cref{thm:existence} under \Cref{ass:tube}. Since $\cD_z$
is open, the path constraint cannot be active at a feasible trajectory. It
therefore has no multiplier and affects neither the representer form nor the
first-order conditions. Problem \eqref{eq:finite} sees $\cD_z$ only through
$\cC_i\subset \cD_z$ at the waypoint times, and Step~1 guarantees validity
between them near a local minimizer. In practice, one solves
\eqref{eq:finite} and checks $\bz^\star(t)\in \cD_z$ along the resulting
trajectory. By the converse of \Cref{thm:reduction}, this a posteriori check
certifies a local minimizer of \eqref{eq:P} without verifying \Cref{ass:tube}. Problem
\eqref{eq:finite} has a convex quadratic objective, constraints on manifolds, and
$nN$ unknowns that do not depend on any time discretization.

\section{Numerical Experiments}\label{sec:num}

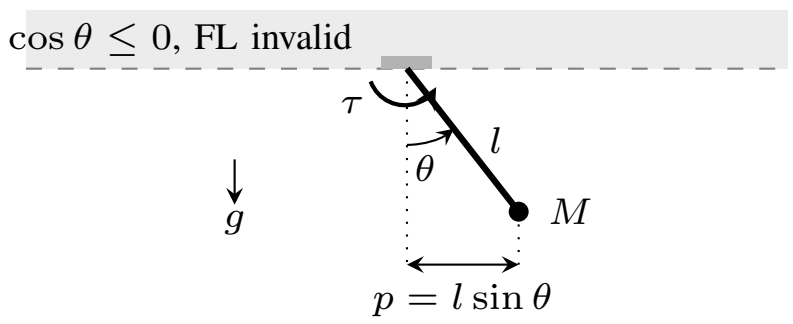
\begin{figure}[!htb]
\centering
\begin{tikzpicture}[scale=1.2,>=stealth]
  \fill[gray!15] (-2.1,0) rectangle (2.1,0.32);
  \node[font=\scriptsize] at (-1.25,0.16) {$\cos\theta\le0$, FL invalid};
  \draw[dashed,gray] (-2.1,0) -- (2.1,0);
  \fill[gray!60] (-0.14,0) rectangle (0.14,0.07);
  \draw[dotted] (0,0) -- (0,-1.12);
  \coordinate (M) at ({sin(38)},{-cos(38)});
  \draw[very thick] (0,0) -- (M) node[midway,right=1pt,font=\scriptsize] {$l$};
  \fill (M) circle (0.055) node[right=2pt,font=\scriptsize] {$M$};
  \draw[->] (0,-0.42) arc[start angle=-90,end angle=-52,radius=0.42];
  \node[font=\scriptsize] at (0.1,-0.55) {$\theta$};
  \draw[->,thick] (-0.2,-0.07) arc[start angle=200,end angle=330,radius=0.2];
  \node[font=\scriptsize] at (-0.3,-0.2) {$\tau$};
  \draw[dotted] (M) -- ({sin(38)},-1.08);
  \draw[<->] (0,-1.08) -- ({sin(38)},-1.08)
    node[midway,below,font=\scriptsize] {$p=l\sin\theta$};
  \node[font=\scriptsize] at (-0.95,-0.85) {$g$};
  \draw[->] (-0.95,-0.5) -- (-0.95,-0.75);
\end{tikzpicture}
\caption{Single-link arm. The end-effector position $p$
linearizes the dynamics only while the arm stays below the horizontal.}
\label{fig:arm}
\end{figure}

We introduce the dynamics of the single-link arm of \Cref{fig:arm} in
Section~\ref{sec:num-model} and the task we solve in
Section~\ref{sec:num-task}. The last two subsections show two results: 1) the
finite-dimensional reduction is exact and has far fewer variables than a fine
discretization; and 2) the linear model $(\mA,\mB)$
can act as a design parameter to keep the optimal trajectory in the region
where feedback linearization is valid and make the minimized cost
approximate the physical one.

\subsection{Nonlinear System and Linear Models}\label{sec:num-model}

A point mass $M$ at the tip of a link of length $l$ is driven by torque $\tau$, with the angle $\theta$ measured from downward,
\begin{equation}\label{eq:arm}
    Ml^2\ddot\theta=\tau-Mgl\sin\theta,\qquad \bx=(\theta,\dot\theta),\quad u=\tau.
\end{equation}
We take the horizontal end-effector position $p=l\sin\theta$ as output.
Then $\bz=\mT(\bx)\coloneqq(l\sin\theta,\,l\cos\theta\,\dot\theta)=(p,\dot p)$ and
\begin{equation}\label{eq:arm-pdd}
\begin{aligned}
    \ddot p&=\chi(\bx)+\gamma(\bx)\,\tau,\qquad \gamma(\bx)\coloneqq\frac{\cos\theta}{Ml},\\
    \chi(\bx)&\coloneqq-l\sin\theta\,\dot\theta^2-g\sin\theta\cos\theta.
\end{aligned}
\end{equation}
The map $\mT$ is a diffeomorphism from $\cD=(-\pi/2,\pi/2)\times\mathbb{R}$ onto
$\cD_z=(-l,l)\times\mathbb{R}$, but not a global one. When the arm is
horizontal, $\gamma$ vanishes and the torque has no effect on $\ddot p$. The
path constraint $|z_1(t)|<l$ thus requires the arm to stay below the
horizontal.


\begin{table}[!htb]
\centering
\caption{The two linear models of the arm.}
\label{tab:models}
\renewcommand{\arraystretch}{2.1}
\renewcommand{\cellset}{\renewcommand{\arraystretch}{1.1}}
\begin{tabular*}{\columnwidth}{@{\extracolsep{\fill}}lcc@{}}
\toprule
 & Brunovsk\'y & Pendulum-matched\\
\midrule
$\Omega$ & $0$ & $\sqrt{g/l}$\\
$\alpha(\bx)$ & $-\dfrac{Ml\,\chi(\bx)}{\cos\theta}$
 & $-\dfrac{Ml\,\big(\chi(\bx)+g\sin\theta\big)}{\cos\theta}$\\
$v$ & $\ddot p$ & $\ddot p+(g/l)\,p$\\
Cost penalizes & \makecell{horizontal\\ acceleration} & \makecell{approximated\\ torque}\\
Spline basis & $1,\ t,\ t^2,\ t^3$
 & \makecell{$\cos\Omega t,\ \sin\Omega t,$\\ $t\cos\Omega t,\ t\sin\Omega t$}\\
\bottomrule
\end{tabular*}
\end{table}

We use a linear model $\mA=\left[\begin{smallmatrix}0&1\\
-\Omega^2&0\end{smallmatrix}\right]$ and $\mB=[0\ 1]\tp$ with the two choices
of $\Omega$ in \Cref{tab:models}. 
The pendulum-matched model follows from the Brunovsk\'y model by the
feedback $\mF=[-g/l\ \ 0]$ of \Cref{prop:fl-AB}. Its virtual input
$v=(\cos\theta/(Ml))\,\tau+g\sin\theta(1-\cos\theta)-l\sin\theta\,\dot\theta^2$
equals $\tau/(Ml)$ to first order around the hanging rest. Its cost
thus matches $(Ml)^{-2}\int_0^{\tfin}\tau^2\,\mathrm{d}t$ up to a relative
error of second order near this rest, whereas the Brunovsk\'y cost penalizes
the net torque $\tau-Mgl\sin\theta$ to first order. The spline bases in \Cref{tab:models}
follow from the eigenvalues of $\mA$ and $-\mA$.

\subsection{Task and Finite-Dimensional Program}\label{sec:num-task}

We solve the following task: the arm starts and ends at rest, $\bx(0)=\bx(\tfin)=\mathbf{0}$, and must reach the joint
speeds $\dot\theta(t_1)=\omega_1$ and $\dot\theta(t_2)=\omega_2$ at angles
left free. 
$\omega_1$ and $\omega_2$ have opposite signs. 
The $i$th of these
waypoints has $\bc_i=[0\ 1]\tp$ and $y_i=\omega_i$. Since $\cos\theta>0$ on
$\cD$, its feasible set
$\cC_i=\{\bzeta\in \cD_z:\ \zeta_2=\omega_i(l^2-\zeta_1^2)^{1/2}\}$ needs only the
positive root and is a nonconvex arc. The final condition $\bx(\tfin)=\mathbf{0}$ consists of two waypoints at
$\tfin$ that prescribe the full state and reduce to $\bz(\tfin)=\mathbf{0}$. With $\bz_0=\mT(\mathbf{0})=\mathbf{0}$, problem \eqref{eq:P} for this task reads
\begin{equation}\label{eq:arm-P}
\begin{aligned}
    \min_{v}\ \ & \int_0^{\tfin}R\,v(t)^2\,\mathrm{d}t\\
    \text{s.t.}\ \ & \dot \bz=\mA\bz+\mB v,\quad \bz(0)=\mathbf{0},\quad |z_1(t)|<l,\\
    & \bz(t_i)\in \cC_i,\ \ i=1,2,\quad \bz(\tfin)=\mathbf{0}.
\end{aligned}
\end{equation}
Waypoints that share a time share one coefficient in \eqref{eq:finite},
since $\mK_1(\cdot,t)\ba+\mK_1(\cdot,t)\ba'=\mK_1(\cdot,t)(\ba+\ba')$. With the distinct
times $\bar t_1=t_1$, $\bar t_2=t_2$, and $\bar t_3=\tfin$, the unknown is
$\ba\in\mathbb{R}^6$, and \eqref{eq:finite} reads
\begin{equation}\label{eq:arm-nlp}
\begin{aligned}
    \min_{\ba\in\mathbb{R}^6}\ \ & \ba\tp\bGamma\ba\\
    \text{s.t.}\ \ & [\bzeta_i(\ba)]_2=\omega_i\big(l^2-[\bzeta_i(\ba)]_1^2\big)^{1/2},
    \ \ i=1,2,\\
    & \bzeta_3(\ba)=\mathbf{0},
\end{aligned}
\end{equation}
where $\bzeta_j(\ba)=\sum_{k=1}^3\mK_1(\bar t_j,\bar t_k)\ba_k$ and $\bGamma$ is built
from the kernel $\mK_1$ of Section~\ref{sec:prelim-math} for each model.

The constraints of \eqref{eq:arm-nlp} read $\bGamma\ba=\bd(\theta_1,\theta_2)$
with $\bd\coloneqq(l\sin\theta_1,\,\omega_1l\cos\theta_1,\,l\sin\theta_2,\,
\omega_2l\cos\theta_2,\,0,\,0)$ and $\theta_i=\theta(t_i)\in(-\pi/2,\pi/2)$.
Hence \eqref{eq:arm-nlp} is equivalent to minimizing
$J(\theta_1,\theta_2)\coloneqq \bd\tp\bGamma\iv \bd$ over two angles. The actual
degrees of freedom are thus the waypoint angles $\theta(t_1)$ and
$\theta(t_2)$. Once they are fixed, the waypoint states are fixed, and the
minimizer is the interpolating spline of Section~\ref{sec:prelim-math}.

We set $M=l=1$, $g=9.81$, $R=1$, $t_1=1.1$, $t_2=2.1$, and $\tfin=3$.
Sections~\ref{sec:num-exact} and~\ref{sec:num-region} use
$(\omega_1,\omega_2)=(1,-3)$. Section~\ref{sec:num-exact} tests the
exactness of the reduction on the pendulum-matched model, whose
trajectories stay in $\cD_z$. Section~\ref{sec:num-design} compares the two
linear models, and its second part uses the slower constraints
$(\omega_1,\omega_2)=(1,-1.5)$.

\subsection{Exactness of the Reduction for the Pendulum-matched Model}\label{sec:num-exact}

We show that the finite-dimensional program \eqref{eq:arm-nlp} and a fine
discretization of \eqref{eq:arm-P} have the same cost landscape and the same
KKT points, and that \eqref{eq:arm-nlp} is much cheaper to solve.

As ground truth, we discretize \eqref{eq:arm-P} with $v$ piecewise constant on a
uniform grid of step $h=\tfin/450$ and exact discrete-time dynamics. We write this
discretized nonlinear programming (NLP) by direct multiple shooting, which is generally preferred for
optimal control problems \cite{bock1984multiple}. Its variables are the states
at the nodes after $t=0$ and the inputs on the $450$ intervals, and the
dynamics enter as equality constraints.
We find the KKT points of \eqref{eq:arm-nlp} and of the discretized
\eqref{eq:arm-P} by Newton's method from the same $16$ initial guesses. For a
guess $(\theta_1^0,\theta_2^0)$, both NLPs start from the interpolating spline
of Section~\ref{sec:prelim-math} through the waypoint states fixed by these
angles, which minimizes $\int_0^{\tfin}v^2\,\mathrm{d}t$. The difference is that $\ba$ parametrizes this
spline exactly, whereas the states and inputs of direct multiple shooting are
a discretized approximation of it. The local curvature classifies each KKT
point as a local minimizer or a saddle point. The path constraint is checked
a posteriori.

\begin{figure}[t]
\centering
\includegraphics[width=\columnwidth]{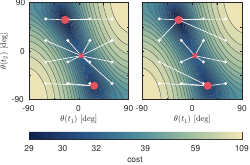}
\caption{Color map of the cost $J$ of \eqref{eq:arm-nlp} (left) and of the
cost $J_h$ of the discretized \eqref{eq:arm-P} (right) over the waypoint
angles $(\theta_1,\theta_2)$ for the pendulum-matched model, on a logarithmic
scale. Dots mark the local minimizers and
crosses the saddle points. White lines join each initial guess to the KKT
point that Newton's method reaches from it.}
\label{fig:arm-exact}
\end{figure}

\begin{figure}[t]
\centering
\includegraphics[width=\columnwidth]{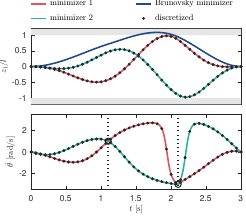}
\caption{End-effector position (upper) and joint speed (lower) for
$(\omega_1,\omega_2)=(1,-3)$. For the pendulum-matched model, the two
minimizers of \eqref{eq:arm-nlp} (lines) and of the discretized
\eqref{eq:arm-P} (dots) coincide, which shows exactness. They stay out of the
shaded invalid region $|z_1|\ge l$, whereas the Brunovsk\'y minimizer enters
it. Circles and dotted lines mark the waypoints $\dot\theta(t_i)=\omega_i$.}
\label{fig:arm-traj}
\end{figure}

\begin{table}[!htb]
\centering
\caption{KKT points of \eqref{eq:arm-nlp} and the discretized
\eqref{eq:arm-P}. $\Delta_{\bz}\coloneqq\|\bz_h-\bz\|/\|\bz\|$
and $\Delta_{\dot\theta}$ are the relative errors between the trajectories
$\bz_h$ of the discretized \eqref{eq:arm-P} and $\bz$ of
\eqref{eq:arm-nlp}.}
\label{tab:arm-kkt}
\renewcommand{\arraystretch}{1.3}
\setlength{\tabcolsep}{2pt}
\small
\begin{tabular*}{\columnwidth}{@{\extracolsep{\fill}}lcccccc@{}}
\toprule
Point & $\theta_1$ & $\theta_2$ & $J$ & $J_h$ & $\Delta_{\bz}$ & $\Delta_{\dot\theta}$\\
\midrule
Min.\ 1 & $-25.0^\circ$ & $58.3^\circ$ & $29.239$ & $29.240$
 & $1.7\times10^{-5}$ & $3.7\times10^{-5}$\\
Min.\ 2 & $27.7^\circ$ & $-64.2^\circ$ & $30.312$ & $30.313$
 & $2.6\times10^{-5}$ & $6.0\times10^{-5}$\\
Saddle & $4.4^\circ$ & $-8.4^\circ$ & $32.217$ & $32.218$
 & $3.5\times10^{-5}$ & $3.9\times10^{-5}$\\
\bottomrule
\end{tabular*}
\end{table}

As \Cref{fig:arm-exact} shows, the two NLPs have the same cost map, with two
local minimizers and one saddle point, which \Cref{tab:arm-kkt} lists.
Minimizer 1 passes $t_1$ on the left of the hanging rest and $t_2$ on the
right, and minimizer 2 does the opposite. The KKT points of the two NLPs
coincide within $0.01^\circ$. Newton's method converges from all $16$ initial
guesses, and different guesses reach different KKT points, as the white lines
in \Cref{fig:arm-exact} show.

As shown by \Cref{fig:arm-traj}, both minimizers stay within $|z_1(t)|<l$ and
meet $\dot\theta(t_i)=\omega_i$. \Cref{tab:arm-kkt} lists the errors
between the trajectories of the two NLPs.

\begin{table}[!htb]
\centering
\caption{Size of each NLP and median cost of Newton's method to reach a
local minimizer, compiled C++ on one core.}
\label{tab:arm-time}
\renewcommand{\arraystretch}{1.3}
\begin{tabular*}{\columnwidth}{@{\extracolsep{\fill}}lcccc@{}}
\toprule
NLP & Variables & Constraints & Newton steps & Time [$\mu$s]\\
\midrule
\eqref{eq:arm-nlp} & $6$ & $4$ & $5$ & $3.7$\\
Discretized \eqref{eq:arm-P} & $1350$ & $904$ & $14$ & $1980$\\
\bottomrule
\end{tabular*}
\end{table}

Beyond being exact, the reduction is far smaller. \Cref{tab:arm-time} compares
the two NLPs. We time both with a C++ implementation of the same Newton
method. The
discretized NLP has three variables and two dynamics constraints per grid
interval, whereas the size of \eqref{eq:arm-nlp} does not depend on any grid.
Newton's method also needs fewer steps on \eqref{eq:arm-nlp}. Although each
Newton step of multiple shooting solves only a banded linear system,
\eqref{eq:arm-nlp} reaches a local minimizer more than $500$ times faster. A finer grid would widen this gap.

\subsection{The Linear Model as a Design Parameter}\label{sec:num-design}

\subsubsection{Staying in the Region of Validity}\label{sec:num-region}

For $(\omega_1,\omega_2)=(1,-3)$ and the Brunovsk\'y model, Newton's method
finds a single local minimizer of \eqref{eq:arm-nlp} at
$(45.6^\circ,70.4^\circ)$ and two saddle points. As shown by the upper panel
of \Cref{fig:arm-traj}, its cubic spline leaves $\cD_z$ at $t=1.52$ and reaches
$\max_t|z_1(t)|=1.09\,l$, so the arm would pass the horizontal. The a
posteriori check of \Cref{thm:reduction} rejects the only minimizer found, so the
reduction certifies no local minimizer of \eqref{eq:arm-P} for this model.
Minimizer 1 of the pendulum-matched model stays within $0.98\,l$. By the
converse of \Cref{thm:reduction}, it is a local minimizer of \eqref{eq:arm-P}.
This example shows that a proper choice of $(\mA, \mB)$ is able to keep the optimal trajectory in the region where
feedback linearization is valid.

\subsubsection{Approximating the Physical Cost}\label{sec:num-cost}

\begin{figure}[t]
\centering
\includegraphics[width=\columnwidth]{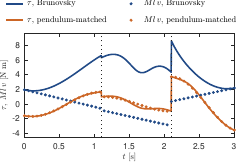}
\caption{Joint torque $\tau$ (lines) and $Ml\,v$ (dots) of the global
minimizers of \eqref{eq:arm-nlp} for $(\omega_1,\omega_2)=(1,-1.5)$. Blue is
the Brunovsk\'y model and orange the pendulum-matched model. Dotted lines
mark $t_1$ and $t_2$. Signals jump there because the input of each kernel
section $\mK_1(\cdot,\bar t_k)\ba_k$ vanishes after $\bar t_k$.}
\label{fig:arm-cost}
\end{figure}

\begin{table}[!htb]
\centering
\caption{Torque and virtual-input costs of the global minimizers. Gap
$=1-(Ml)^2\!\int v^2\,\mathrm{d}t/\!\int\tau^2\,\mathrm{d}t$.}
\label{tab:arm-cost}
\renewcommand{\arraystretch}{1.4}
\begin{tabular}{@{}lccc@{}}
\toprule
Model & $\int_0^{\tfin}\tau^2\,\mathrm{d}t$ & $(Ml)^2\int_0^{\tfin}v^2\,\mathrm{d}t$ & Gap\\
\midrule
Brunovsk\'y & $64.42$ & $6.80$ & $89.45\%$\\
Pendulum-matched & $7.86$ & $7.78$ & $1.02\%$\\
\bottomrule
\end{tabular}
\end{table}

For $(\omega_1,\omega_2)=(1,-1.5)$, the global minimizers of
\eqref{eq:arm-nlp} stay within $0.89\,l$ for the Brunovsk\'y model and
$0.41\,l$ for the pendulum-matched model. Both are thus certified by the
converse of \Cref{thm:reduction}.
\Cref{fig:arm-cost} compares their torques with the virtual inputs, and
\Cref{tab:arm-cost} compares the costs. For the pendulum-matched model,
$Ml\,v$ follows $\tau$ closely, and the virtual-input cost is within $1.0\%$ of
the torque cost. For the Brunovsk\'y model, the virtual-input cost misses
$89.5\%$ of the torque cost. Its cost penalizes the net torque
$\tau-Mgl\sin\theta$ rather than $\tau$, so the planner does not see the
torque that holds the arm against gravity. Its minimizer swings up to
$63^\circ$ and needs eight times the integrated squared torque of the pendulum-matched
minimizer.
This example shows that a proper choice of $(\mA, \mB)$ is able to make the minimized cost approximate the
physical one better.

\section{Conclusion}\label{sec:conclusion}

We extended the exact finite-dimensional reduction of WCTO from linear
systems to full-state feedback-linearizable systems with waypoints imposed on
the physical state. Every local minimizer is a finite kernel expansion. A
finite-dimensional program with $nN$ unknowns recovers it, and an a
posteriori check certifies its solutions. Experiments on a single-link arm
confirmed that the reduction is exact. They also showed that the linear model
$(\mA,\mB)$ offers degrees of freedom. On the arm, a choice guided by domain
knowledge kept the solution where the linearization is valid and made the cost
approximate the integrated squared torque. Future
work will address the systematic choice of $(\mA,\mB)$.

\section*{Acknowledgment}

Claude (Anthropic) helped draft and edit this manuscript under the direction
and review of the authors, who are responsible for the results.

\bibliographystyle{IEEEtran}
\bibliography{refs}

@book{khalil2002nonlinear,
  author    = {Khalil, Hassan K.},
  title     = {Nonlinear Systems},
  edition   = {3rd},
  publisher = {Prentice Hall},
  address   = {Upper Saddle River, NJ},
  year      = {2002}
}

@article{zhang1997splines,
  author  = {Zhang, Z. and Tomlinson, J. and Martin, C.},
  title   = {Splines and linear control theory},
  journal = {Acta Applicandae Mathematicae},
  volume  = {49},
  number  = {1},
  pages   = {1--34},
  year    = {1997}
}

@book{egerstedt2010control,
  author    = {Egerstedt, Magnus and Martin, Clyde},
  title     = {Control Theoretic Splines: Optimal Control, Statistics, and Path Planning},
  publisher = {Princeton University Press},
  address   = {Princeton, NJ},
  year      = {2010}
}

@article{aubin2021linearly,
  author  = {Aubin-Frankowski, Pierre-Cyril},
  title   = {Linearly constrained linear quadratic regulator from the viewpoint of kernel methods},
  journal = {SIAM Journal on Control and Optimization},
  volume  = {59},
  number  = {4},
  pages   = {2693--2716},
  year    = {2021}
}

@article{guillemet2025adaptive,
  title = {Adaptive Vector-Valued Splines for the Resolution of Inverse Problems},
  author = {Guillemet, Vincent and Unser, Michael},
  year = 2027,
  month = jan,
  journal = {Journal of Computational and Applied Mathematics},
  volume = {489},
  pages = {117910},
  issn = {0377-0427},
  doi = {10.1016/j.cam.2026.117910},
  urldate = {2026-07-01}
}

@inproceedings{mellinger2011minimum,
  author    = {Mellinger, Daniel and Kumar, Vijay},
  title     = {Minimum snap trajectory generation and control for quadrotors},
  booktitle = {Proc. IEEE Int. Conf. Robot. Autom.},
  pages     = {2520--2525},
  year      = {2011}
}

@article{wang2022geometrically,
  author  = {Wang, Zhepei and Zhou, Xin and Xu, Chao and Gao, Fei},
  title   = {Geometrically constrained trajectory optimization for multicopters},
  journal = {IEEE Transactions on Robotics},
  volume  = {38},
  number  = {5},
  pages   = {3259--3278},
  year    = {2022}
}

@article{brunovsky1970classification,
  author  = {Brunovsk{\'y}, Pavol},
  title   = {A classification of linear controllable systems},
  journal = {Kybernetika},
  volume  = {6},
  number  = {3},
  pages   = {173--188},
  year    = {1970}
}

@article{kalman1963controllability,
  author  = {Kalman, Rudolf E. and Ho, Yu-Chi and Narendra, Kumpati S.},
  title   = {Controllability of linear dynamical systems},
  journal = {Contributions to Differential Equations},
  volume  = {1},
  number  = {2},
  pages   = {189--213},
  year    = {1963}
}

@article{flash1985coordination,
  author  = {Flash, Tamar and Hogan, Neville},
  title   = {The coordination of arm movements: An experimentally confirmed mathematical model},
  journal = {Journal of Neuroscience},
  volume  = {5},
  number  = {7},
  pages   = {1688--1703},
  year    = {1985}
}

@article{lin1983formulation,
  author  = {Lin, Chun-Shin and Chang, Po-Rong and Luh, J. Y. S.},
  title   = {Formulation and optimization of cubic polynomial joint trajectories for industrial robots},
  journal = {IEEE Transactions on Automatic Control},
  volume  = {28},
  number  = {12},
  pages   = {1066--1074},
  year    = {1983}
}

@inproceedings{werling2010optimal,
  author    = {Werling, Moritz and Ziegler, Julius and Kammel, S{\"o}ren and Thrun, Sebastian},
  title     = {Optimal trajectory generation for dynamic street scenarios in a {F}ren\'et frame},
  booktitle = {Proc. IEEE Int. Conf. Robot. Autom.},
  pages     = {987--993},
  year      = {2010}
}

@article{bock1984multiple,
  author  = {Bock, H. G. and Plitt, K. J.},
  title   = {A multiple shooting algorithm for direct solution of optimal control problems},
  journal = {IFAC Proceedings Volumes},
  volume  = {17},
  number  = {2},
  pages   = {1603--1608},
  year    = {1984}
}

@article{sidhu1979vector,
  author  = {Sidhu, G. S. and Weinert, H. L.},
  title   = {Vector-valued {L}g-splines {I}. {I}nterpolating splines},
  journal = {Journal of Mathematical Analysis and Applications},
  volume  = {70},
  number  = {2},
  pages   = {505--529},
  year    = {1979}
}

@article{sidhu1984vector,
  author  = {Sidhu, G. S. and Weinert, H. L.},
  title   = {Vector-valued {L}g-splines {II}. {S}moothing splines},
  journal = {Journal of Mathematical Analysis and Applications},
  year    = {1984},
  doi     = {10.1016/0022-247X(84)90109-4}
}

@book{biagiotti2008trajectory,
  author    = {Biagiotti, Luigi and Melchiorri, Claudio},
  title     = {Trajectory Planning for Automatic Machines and Robots},
  publisher = {Springer},
  address   = {Berlin, Heidelberg},
  year      = {2008}
}

@article{beaver2024optimal,
  author  = {Beaver, Logan E. and Malikopoulos, Andreas A.},
  title   = {Optimal control of differentially flat systems is surprisingly easy},
  journal = {Automatica},
  volume  = {159},
  pages   = {111404},
  year    = {2024}
}

@article{unser2017splines,
  title={Splines are universal solutions of linear inverse problems with generalized 
  {TV} regularization},
  author={Unser, Michael and Fageot, Julien and Ward, John Paul},
  journal={SIAM Review},
  volume={59},
  number={4},
  pages={769--793},
  year={2017},
  publisher={SIAM}
}

@book{wahba1990spline,
  title={Spline models for observational data},
  author={Wahba, Grace},
  year={1990},
  publisher={SIAM}}

@inproceedings{scholkopf2001generalized,
  title={A generalized representer theorem},
  author={Sch{\"o}lkopf, Bernhard and Herbrich, Ralf and Smola, Alex J},
  booktitle={International conference on computational learning theory},
  pages={416--426},
  year={2001},
  organization={Springer}
}

@article{de1966splines,
  title={On splines and their minimum properties},
  author={De Boor, Carl and Lynch, Robert E},
  journal={Journal of Mathematics and Mechanics},
  volume={15},
  number={6},
  pages={953--969},
  year={1966},
  publisher={JSTOR}}

@article{schultz1967splines,
  title={L-splines},
  author={Schultz, Martin H and Varga, Richard S},
  journal={Numerische Mathematik},
  volume={10},
  number={4},
  pages={345--369},
  year={1967},
  publisher={Springer}
}

@article{unser2005cardinal,
  title={Cardinal exponential splines: Part {I}-{T}heory and filtering algorithms},
  author={Unser, Michael and Blu, Thierry},
  journal={IEEE Transactions on Signal Processing},
  volume={53},
  number={4},
  pages={1425--1438},
  year={2005},
  publisher={IEEE}
}

\end{document}